\documentclass[11pt,twoside]{article}
\usepackage{enumerate}
\usepackage{amsmath}
\usepackage{amssymb}
\usepackage{amsthm}
\usepackage{amscd}
\usepackage{catmac}
\DeclareSymbolFontAlphabet{\Bbb}{AMSb}
\newcommand{\id}{\text{{\rm id}}}

\newcommand{\codim}{\text{{\rm codim}}}

\newcommand{\norm}[1]{ \|#1 \| }
\newcommand{\bignorm}[1]{\left \|#1 \right \| }
\newcommand{\convex}[1]{\mathbf{M^{(#1)}}}
\newcommand{\concave}[1]{\mathbf{M_{(#1)}}}

\newcommand{\K}{\Bbb{K}}

\newcommand{\R}{\Bbb{R}}

\newcommand{\N}{\Bbb{N}}

\newcommand{\LL}{{\cal L}}

\newcommand{\Id}{\hookrightarrow}

\newcommand{\ui}{\mathcal{S}}
\theoremstyle{definition}
\newtheorem{defin}{Definition}[section]

\theoremstyle{plain}

\newtheorem{lemma}[defin]{Lemma}
\newtheorem{theo}[defin]{Theorem}
\newtheorem{cor}[defin]{Corollary}
\newtheorem{prop}[defin]{Proposition}

\theoremstyle{remark}
\numberwithin{equation}{section}
\begin{document}
\bibliographystyle{amsalpha}
\title{\bf Summing inclusion maps between symmetric sequence spaces}
\author{by \\[20pt] Andreas Defant, Mieczys{\l}aw Masty{\l}o\thanks{Research supported by KBN Grant 2 P03A
042 18}\,\, and Carsten Michels \\[10pt]  }
\date{}
  \maketitle
\begin{abstract}
In 1973/74 Bennett and (independently) Carl proved that for $1 \le u \le 2$ the identity map
id: $\ell_u \Id
\ell_2$ is absolutely $(u,1)$-summing, i.\,e. for every unconditionally summable sequence $(x_n)$
in $\ell_u$ the
scalar sequence $(\|x_n \|_{\ell_2})$ is contained in $\ell_u$,
 which improved upon well-known results of Littlewood and Orlicz. The following substantial
 extension is our main result:
 For a $2$-concave
symmetric Banach sequence space $E$ the identity map $\text{id}: E \Id \ell_2$ is absolutely
$(E,1)$-summing, i.\,e. for every unconditionally summable sequence $(x_n)$ in $E$ the scalar sequence
$(\|x_n \|_{\ell_2})$ is contained in
$E$. Various applications are given, e.\,g. to the theory of eigenvalue distribution of compact
operators where we show
that the sequence of
eigenvalues of an operator $T$ on $\ell_2$ with values in a $2$-concave symmetric Banach
sequence space $E$ is a
multiplier from $\ell_2$ into $E$. Furthermore,
 we prove an asymptotic formula for the $k$-th approximation number
 of the identity map $\text{id}: \ell_2^n \Id E_n$, where $E_n$ denotes the linear span of the
 first $n$ standard unit vectors in $E$, and apply it to
 Lorentz  and Orlicz sequence spaces.
\end{abstract}
\section{Introduction}
In 1930 Littlewood \cite{little} proved that for every bilinear and continuous operator
\mbox{$\varphi: c_0 \times c_0 \rightarrow \R$} the quantity
$\sum_{k,\ell=1}^\infty |\varphi(e_k,e_\ell)|^{4/3}$ is finite; this is
equivalent to the statement that for every unconditionally summable sequence $(x_n)$ in $\ell_1$
the scalar sequence $(\norm{x_n}_{\ell_{4/3}})$ is contained in $\ell_{4/3}$. Bennett \cite{bennett}
and (independently) Carl \cite{carl} extended Littlewood's result in the following way: For
$1 \le u \le v \le 2$ and every unconditionally summable sequence $(x_n)$ in $\ell_u$ the
sequence $(\norm{x_n}_{\ell_v})$ is contained in $\ell_r$, where $1/r =1/u-1/v+1/2$. Their result
has useful applications in various parts of analysis---in particular, in approximation theory
as well as for the theory of eigenvalue distribution of
compact operators, e.\,g. that for $1\le u <2$ every operator on $\ell_2$ with values in $\ell_u$
has absolutely $r$-summable eigenvalues, where $1/r=1/u-1/2$.
\par
The case $v=2$ in the Bennett--Carl result is crucial (for the proof as well as for applications).
 Motivated by applications
to interpolation theory (see e.\,g. \cite{ovchinni} and \cite{milman}) Maligranda and the
second named author in \cite{mastylo} proved that for an Orlicz function $\varphi$ for which the
map $t \mapsto \varphi(\sqrt{t})$ is equivalent to a concave function and for every unconditionally
summable sequence $(x_n)$ in the Orlicz sequence space $\ell_\varphi$ the sequence
$(\norm{x_n}_{\ell_2})$ is contained in $\ell_\varphi$. Moreover, based on complex interpolation, in
\cite{dm98} various commutative and non commutative variants were given.
\par
These results were the starting point for
the research on which this article is based upon. Developing and using
complex interpolation formulas for spaces of operators related to those of Kouba \cite{kouba},
our main result is a far reaching extension of the
 above results: For a $2$-concave
 symmetric Banach sequence space $E$  and every unconditionally summable sequence $(x_n)$ in
 $E$ the sequence $(\norm{x_n}_{\ell_2})$ is contained in $E$. In the language of
 $(E,1)$-summing operators (which we will recall later on) this means that the identity map
 \mbox{$\id : E \Id \ell_2$} is  $(E,1)$-summing. An example shows that the $2$-concavity of
 $E$ is not superfluous.  As in the classical case our result
  has some useful applications.  We show that the sequence of eigenvalues
   of an operator $T$ on $\ell_2$ with values in a $2$-concave symmetric Banach sequence space
   $E$ is a multiplier from $\ell_2$ into $E$, a result which for $E=\ell_u$, $1 \le u \le 2$, is
   well-known (note that the space of multipliers from $\ell_2$ into $\ell_u$ coincides with
 $\ell_r$, $1/r =1/u-1/2$). Furthermore, we prove for a $2$-concave
   symmetric Banach sequence space $E$ and
   $1 \le k \le n$ the asymptotic formula
$$
a_k(\id: \ell_2^n \Id E_n) \asymp \frac{\lambda_E(n-k+1)}{(n-k+1)^{1/2}},
$$
where $a_k(T)$ denotes the $k$-th approximation number of an operator $T$, $E_n$ stands for
the linear span of the first $n$ standard unit vectors in $E$ and $\lambda_E: \N \rightarrow
\R_+$ is the fundamental function of the sequence space $E$, and apply it to Lorentz and
Orlicz sequence spaces.
\section{Preliminaries}
For a positive number $a$ we denote by $\lfloor a \rfloor$ the largest integer less or equal than $a$.
If $(a_n)$ and $(b_n)$ are scalar sequences we write \mbox{$a_n \prec b_n$}
 whenever there is some $c \ge 0$ such that $a_n \le c \cdot b_n$ for all $n$,
  and
\mbox{$a_n \asymp b_n$} whenever \mbox{$a_n \prec b_n$} and
\mbox{$b_n \prec a_n$}. \par
 We  use standard notation and
 notions from Banach space theory, as presented e.\,g. in
 \cite{lt77}, \cite{lt} and \cite{tj}.
If $E$ is a Banach space, then
 $B_E$ is its (closed) unit ball and $E'$ its dual space.
\par
 Throughout the paper by a Banach sequence
space we mean a real Banach lattice $E$ modelled on the set of positive
integers $\N$ which contains an element $x$ with $\mbox{supp} \,x= \N$.
A Banach sequence space $E$ is said to be symmetric provided that
$\|(x_n)\|_{E} = \|(x_n^*)\|_{E}$, where $(x_n^*)$ denotes the decreasing
rearrangement of the sequence $(x_n)$, i.e.
$$
x_{n}^{*}:= \inf \{\sup_{i\in \N \setminus J} |x_i| \, | \,
J\subset \N, \, \mbox{card} (J) < n \}.
$$
 It is maximal if the unit ball $B_E$ is closed in the
pointwise convergence topology induced by the space $\omega$ of all real sequences. Note that this condition is
equivalent to $E^\times = E'$, where as usual $$ E^{\times}:= \{x=(x_n) \in \omega \, |\,  \Sigma_{n=1}^{\infty} |x_n y_n| <
\infty \, \mbox{for all} \, y=(y_n) \in E\} $$ is the K\"othe dual of $E$. Note that $E^{\times}$ is a maximal (symmetric,
provided that $E$ is) Banach sequence ´space under the norm $$ \|x\|:= \sup\{\Sigma_{n=1}^{\infty}|x_n y_n| \, | \,
\|y\|_{E} \leq 1\}. $$ The fundamental function of a symmetric Banach sequence space $E$ is defined by $$
\lambda_{E}(n):=\norm{\textstyle{\sum_{i=1}^{n}} e_i}_E, \quad n \in \N; $$ throughout the paper $(e_n)$ will denote
the standard unit vector basis in $c_0$ and $E_n$ the linear span of the first $n$ unit vectors.
 It is well-known that any symmetric Banach
sequence space $E$ is continuously embedded in the symmetric Marcinkiewicz sequence
space $m_{\lambda_E}$ of all sequences $x=(x_n)$ such that
$$
\|x\|_{\lambda_{E}}:=\sup_{n\geq 1} x_{n}^{**} \lambda_E(n)  < \infty,
$$
where $x_{n}^{**}:=\frac{1}{n} \sum_{k=1}^{n}\, x_{k}^{*}$.
 For the notions of $p$-convexity and $q$-concavity ($1 \le p,q \le \infty$) of a Banach lattice
 $X$ (the associated constants are denoted by $\convex{p}(X)$ and $\concave{q}(X)$, respectively) we refer
 to \cite[1.d.3]{lt}---but since the notion of $2$-concavity is crucial for our purposes recall
 that a Banach sequence space $E$ is called $2$-concave if there exists a constant
 $C>0$ such that for all
 $x_1, \ldots ,x_n \in E$
$$
 \left ( \sum_{i=1}^n \norm{x_i}_E^2 \right )^{1/2}
  \le C \cdot \bignorm{\left(\sum_{i=1}^n |x_i|^2 \right )^{1/2}}_E.
$$
It is well-known that this is equivalent to the notion of cotype~$2$ (see \cite[1.f.16]{lt});
recall that a Banach space $X$ has cotype $q$
($2\leq q <\infty $) if there is a constant $C>0$ such that for
finitely many $x_1, \ldots, x_n \in X$
\begin{eqnarray*}
\left ( \sum_{i=1}^n \norm{x_i}_X^q \right )^{1/q}
\le C \cdot \Big (\int_{0}^{1} \big \|\sum_{i=1}^n r_{i}(t) \cdot x_i
\big \|_{X}^{2}\, dt \Big )^{1/2}.
\end{eqnarray*}
Note that  $2$-concave symmetric Banach sequence spaces are separable and
maximal.
  An important tool for our purposes are powers of sequence spaces: Let $E$ be a (maximal) symmetric Banach sequence space
  and $0<r<\infty$ such that $\convex{\max(1,r)}(E)=1$. Then
$$
E^r := \{ x \in \ell_\infty \, | \, |x|^{1/r} \in E \}
$$
endowed with the norm
$$
\norm{x}_{E^r}:= \norm{|x|^{1/r}}_E^r, \quad x \in E^r
$$
is again a (maximal) symmetric Banach sequence space which is $1/\min(1,r)$-convex.
For two Banach sequence spaces $E$ and $F$ the space of multipliers $M(E,F)$
from $E$ into $F$ consists of all scalar sequences $x=(x_n) $ such that
the associated multiplication operator $(y_n) \mapsto (x_n\,y_n)$ is
defined and bounded from $E$ into $F$. $M(E,F)$ is a (maximal symmetric
provided that $E$ and $F$ are) Banach sequence space equipped with the norm
$$
\|x\|_{M(E,F)} :=\sup\{\|xy\|_F \, | \, y \in B_{E}\}.
$$
Note that if $E$ is a Banach sequence space then $M(E,\ell_{1})=E^\times$.
 In the case where $E=\ell_2$ and $F$ is
 $2$-concave with $\concave{2}(F)=1$ it can be
 easily seen that
\begin{equation}
\label{diagonalpower}
M(\ell_2,F) = (((F^\times)^2)^\times)^{1/2}
\end{equation}
holds isometrically. We will need that for any symmetric Banach sequence space $E \Id \ell_2$ not equivalent to
$\ell_2$
\begin{equation}
\label{czero}
M(\ell_2,E) \hookrightarrow c_0.
\end{equation}
In fact, for $F:=M(\ell_2,E)$, by the assumption we have
$$
 \lim_{n \rightarrow \infty}
 \lambda_F(n)=\sup_n \norm{\textstyle{\sum_1^n e_i}}_F =\sup_n \norm{\id : \ell_2^n \Id E_n}=
  \infty.
$$
Since for any $x=(x_n) \in F$ the estimate
$
x_n^*\cdot \lambda_F(n) \le \norm{x}_F
$
holds, the claim follows. 
\par For all information on Banach operator ideals and $s$-numbers see \cite{djt}, \cite{koenig},
\cite{pietsch} and \cite{eigenvalues}.
As usual $\LL(E,F)$
 denotes the Banach space of all (bounded and linear) operators from
 $E$ into $F$ endowed with the operator norm $\norm{\cdot}$.
 For an operator $T: X \rightarrow Y$ between Banach spaces recall the definition of
the $k$-th approximation number
$$ a_k(T):= \inf \{ \norm{T-T_k} \, | \, T_k \in \LL(X,Y) \text{ has rank} \,<k\},$$
the $k$-th Weyl number
$$
x_k(T):= \sup \{a_k(TS) \, | \, S \in \LL(\ell_2,X) \text{ with } \norm{S}\le 1\}
$$
and the $k$-th Gelfand number
$$
c_k(T) := \inf \{\norm{T_{|G}} \, | \, G \subset X, \, \codim\, G <k\}.
$$
\par
Moreover, for an $s$-number function
 $s$ and a maximal symmetric Banach sequence space $E$ we denote by $\ui_E^s$ the Banach operator ideal
  of all operators $T$ with $(s_n(T)) \in E$,
 endowed with the norm $\norm{T}_{\ui_E^s}:= \norm{(s_n(T))}_E$; on $\ell_2$ and for fixed $E$ all
 these ideals coincide (isometrically)---for simplicity we then
  denote this space by $\ui_E$.
\par For basic results and notation from interpolation theory we refer to
\cite{bk} and \cite{BL}. We recall that a mapping $\cal F$ from (a subclass $\mathcal{C}$ of) the category of all
couples of Banach spaces into the category
 of all Banach spaces is said to be a method of interpolation (on $\mathcal{C}$)
if for any couple $(X_0,X_1)$ ($\in \mathcal{C}$), the Banach space
${\cal F}(X_0, X_1)$ is intermediate
with respect to $(X_0, X_1)$
(i.\,e. \mbox{$X_0 \cap X_1 \Id {\cal F}(X_0, X_1) \Id
X_0 +X_1$),} and  $T: {\cal F}(X_0, X_1) \to {\cal F}(Y_0, Y_1)$ for all
Banach couples $(X_0, X_1)$, $(Y_0, Y_1)$ ($\in \mathcal{C}$) and every $T: (X_0, X_1) \to
(Y_0, Y_1)$. Here as usual the notation $T: (X_0, X_1) \to (Y_0, Y_1)$
means that $T: X_0 +X_1 \to Y_0 + Y_1 $ is a linear operator such
that  for $j=0,1$ the restriction of $T$ to the space $X_{j}$ is a bounded operator from
$X_j$ into $Y_{j}$. If additionally
$$
\|T: {\cal F}(X_0, X_1) \to {\cal F}(Y_0, Y_1)\| \leq
\max \{ \|T :X_0 \to Y_0 \|, \, \|T: X_1 \to Y_1 \| \}
$$
holds, then ${\cal F}$ is called an exact method of interpolation (on $\mathcal{C}$).
Concrete examples of exact
interpolation methods are the real method of interpolation $(\cdot,\cdot)_{\theta,p}$,
$0<\theta<1$, $1 \le p \le \infty$ (see
e.\,g. \cite[Chapter~3]{BL}) defined on the class of all Banach couples and
the complex method of interpolation $[\cdot,\cdot]_\theta$, $0<\theta<1$
(see e.\,g. \cite[Chapter~4]{BL}) defined on the class of couples of complex Banach spaces.
Both methods are of power type $\theta$, i.\,e. if $\mathcal{F}=(\cdot,\cdot)_{\theta,p}$ or
$\mathcal{F}=[\cdot,\cdot]_\theta$ then for all
$T:(X_0,X_1) \rightarrow (Y_0,Y_1)$ it holds
\begin{equation}
\label{powerreal}
\norm{T:\mathcal{F}(X_0,X_1) \rightarrow \mathcal{F}(Y_0,Y_1)} \le
\norm{T:X_0 \to Y_0}^{1-\theta}
 \cdot \norm{T:X_1 \to Y_1}^\theta.
\end{equation}
In order to avoid misunderstandings,  if we
interpolate between real Banach spaces using the complex method of
interpolation we mean that we use any interpolation functor which
is an extension of the complex method. For such a functor we use the
original notation $[\cdot, \cdot]_{\theta}$.
\par In what follows we will often use the following well-known fact
(see e.\,g. \cite[2.5.1]{BL}) that for any interpolation space
$X$ with respect to $(X_0, X_1)$ there exists an exact interpolation
functor ${\cal F}$ such that ${\cal F}(X_0, X_1) =X$ up to equivalent
norms. An important class of interpolation spaces are $K$-spaces. Recall
that an intermediate Banach space $X$ with respect to a couple $(X_0, X_1)$
is called a relative $K$-space if, whenever $x\in X$ and $y\in X_0 + X_1$
satisfy
$$
K(t, y; X_0, X_1) \leq K(t, x; X_0, X_1) \, \, \, \mbox{for all } t>0,
$$
then it follows that $y \in X$, where
$$
K(t, x; X_0,X_1):= \inf \, \{ \|x_0 \|_{X_0} + t \|x_1 \|_{X_1} \,|
\, x = x_0 + x_1 \}, \quad  t >0
$$
is the Peetre $K$-functional.
\par A Banach couple $(X_0, X_1)$ is said to be a relative Calder\'on
couple if all interpolation spaces with respect to $(X_0, X_1)$ are
also relative $K$-spaces. This is equivalent to: For each pair of elements
$x \in X_0 + X_1$ and $y\in X_0 + X_1$ satisfying
$K(t,y; X_0, X_1) \leq K(t, x; X_0, X_1)$ for all $t>0$, there exists an
operator $T: (X_0, X_1) \to (X_0, X_1)$ such that $Tx=y$.
\section{$\pmb{(E,p)}$-summing operators}
The following definition is a natural extension of the notion of absolutely $(r,p)$-summing
 operators. For two Banach spaces $E$ and
 $F$ we mean by
 $E \Id F$ that $E$ is contained in $F$, and the natural identity map is continuous; in this case
 we put \mbox{$c_E^F:= \norm{\id: E \Id F}$}
 and $c_p^F := c_{\ell_p}^F$ whenever $\ell_p \Id F$. If $E$ and $F$ are Banach sequence spaces
 with  $\norm{e_n}_E=1$ for all $n$, then obviously
 $\ell_1 \Id E$ and $c_1^E=1$. Note also that $E^\times \Id M(E,F)$ with $c_{E^\times}^{M(E,F)}=1$.
\begin{defin}
\label{Epsumming} For $1 \le p < \infty$ let $E$ be a Banach sequence space such that
$\ell_p \Id E$ and $\norm{e_n}_E=1$ for all $n$.
Then an operator $T:X \rightarrow Y$ between
 Banach spaces $X$ and $Y$ is called $(E,p)$-summing (shortly: $T \in
\Pi_{E,p}$) if there exists a constant $C>0$ such that for all
$x_1, \ldots, x_n \in X$
$$
\norm{(\norm{Tx_i}_Y)_{i=1}^n}_E \le C \cdot c_p^E \cdot \sup_{x' \in B_{X'}} \left
( \sum_{i=1}^n |\langle x',x_i \rangle|^p \right )^{1/p},
$$
where in the sequel $(\xi_i)_{i=1}^n$ denotes the sequence $\sum_{i=1}^n \xi_i \cdot e_i$.
We write $\pi_{E,p}(T)$ for the smallest constant $C$ with the above property; in this way we
 obtain the
 Banach operator ideal $(\Pi_{E,p}, \pi_{E,p})$
(see also \cite{milman}), and for $E=\ell_r$ ($r\ge p$)
the well-known Banach operator ideal $(\Pi_{r,p},\pi_{r,p})$ of all absolutely $(r,p)$-summing
 operators.
\end{defin}
\par Let us collect some later needed observations which are all modelled along classical results
on $(r,p)$-summing operators.
We start with the following simple fact that for each maximal Banach
sequence space $E$ an operator  $T: X \to Y $ is $(E,p)$-summing if and only if the
induced linear operator
$$
\widehat{T}: \ell_p^w(X) \rightarrow E(Y), \quad \widehat{T}(x_n):=(Tx_n)
$$
is defined (and hence bounded). In this case,
$\|\widehat{T}: \ell_{p}^{w}(X) \to E(Y)\|=\pi_{E,p}(T)$ provided that
$c_{p}^{E}=1$. Here and in what follows for a given Banach space $X$,
$\ell_{p}^{w}(X)$ and $E(X)$ denotes the Banach space of all weakly
$p$-summable and absolutely $E$-summable sequences $x=(x_n)$ in $X$
equipped with the norms
$$
\|x\|_{\ell_{p}^{w}(X)}:=
\sup_{x'\in B_{X'}} \Big ( \sum_{n=1}^{\infty}
|\langle x', x_n \rangle|^{p}\Big )^{1/p}
$$
and
$$
\|x\|_{E(X)}:= \|(\|x_n\|_{X})\|_{E},
$$
respectively. It is well-known that the Pietsch Domination Theorem implies that any
\mbox{$p$-summing} operator $T: X \to Y$, $1 \le p <\infty$ is a Dunford-Pettis operator, i.\,e.
$T$ transforms weakly convergent sequences into norm convergent sequences,
and thus by Rosenthal's $\ell_1$-Theorem it is compact whenever $X$ does not contain
a copy of $\ell_{1}$. In general this is not true for $(r,p)$-summing
operators as has been noted by Bennett \cite{bennett}, namely the inclusion
map $\ell_{r} \Id \ell_{\infty}$ is $(r,1)$-summing for any $1<r<\infty$,
however not compact. But even in our  more general case the situation becomes
more favorable for operators acting between special Banach spaces (see also Corollary~\ref{calkin}).
\begin{samepage}
\begin{lemma}
\label{compactm} Let $Y$ be a Banach space and $E$  a Banach
sequence space with $\|e_n\|_E=1$ for all $n$. Then the following holds true:
\begin{enumerate}
\item [(a)] If $T \in \Pi_{E,p }(\ell_{p'}, Y)$ with $1<p<\infty$
and $\ell_p \Id E \Id c_{0}$, then
$T$ is a compact operator.
\item [(b)] If $T \in \Pi_{E, 1}(c_0, Y)$ with $\ell_1 \Id
E \Id c_0$, then $T$ is a compact operator.
\end{enumerate}
\end{lemma}
\end{samepage}
\proof
(a) Suppose $T$ is not compact. Then $T$ is no Dunford-Pettis operator
by the reflexivity of $\ell_{p'}$. Thus there exists a sequence
$(x_n)$ in $ \ell_{p'}$ such that $x_n \to 0$ weakly
and $\|Tx_n\|_{Y} \geq C$ for all $n$ with some constant $C>0$.
In consequence $Tx_n \to 0$ weakly in $Y$ and $\|x_n\|_{\ell_{p'}} \geq C/\|T\|$.
Passing to a subsequence, we may assume by the  Bessaga-Pe{\l}czy\'nski Selection
Theorem that $(x_{n})$ is equivalent to a block basis of the unit vector
basis in $\ell_{p'}$ and thus to the unit vector basis in $\ell_{p'}$.
Then $(x_n)$ is weakly $p$-summable in $\ell_{p'}$ since clearly $(e_n)$ is.
But $T: \ell_{p'} \to Y$ is $(E,p)$-summing, hence $(\|Tx_n\|_Y) \in E$, and in particular
$(\|Tx_n\|_Y) \in c_0$, a contradiction.
For (b) we similarly show that $T: c_0 \to Y$ is a Dunford-Pettis operator,
and thus compact since $c_0$ does not contain a copy of $\ell_1$. \qed
\par
 In the following three lemmas we fix $1 \le p < \infty$, and $E$
will always be a Banach sequence space such that $\ell_p \Id E$ and $\norm{e_n}_E=1$ for all $n$.
\begin{lemma}
\label{composition} For an operator $T:X \rightarrow Y$
 between Banach spaces the following are equivalent:
\begin{enumerate}[(a)]
\item $T \in \Pi_{E,p}$, and $\pi_{E,p}(T) \le C$.
\item For all $m$ the map $\Phi^m(T): \LL(\ell_2^m,X) \rightarrow E_m(Y), \quad S \mapsto (TSe_i)$ has norm
$\le C$.
\item $\pi_{E,p}(TS) \le C$ for all $m$ and $S \in \LL(\ell_{p'}^m,X) $ with $\norm{S} \le 1$.
\end{enumerate}
In particular, in this case,
\begin{equation}
\label{discrete}
\pi_{E,p}(T) = \sup_m \norm{\Phi^m(T)} = \sup_m \{ \pi_{E,p}(TS) \, | \, \norm{S: \ell_2^m
\rightarrow X} \le 1 \}.
\end{equation}
\em
The proof follows immediately from the definition and the standard observation
that for each \mbox{$S=\sum_{j=1}^m e_j \otimes x_j \in \LL(\ell_{p'}^m, X)$}
$$
\norm{S} = \sup_{x' \in B_{X'}} \left ( \sum_{j=1}^m |\langle x',x_j \rangle|^p \right
)^{1/p}.
$$
\end{lemma}
\par The following is an analogue of the well-known inclusion formulas (in the classical case due
to Kwapie\'{n} \cite{kwapien} and Tomczak-Jaegermann \cite{tj70}).
\begin{lemma}
\label{inclusion} For $1 \le p <q<\infty$ let $1 < r<\infty$ such that $1/r=1/p-1/q$. Then
 $$ \Pi_{E,p} \subset \Pi_{M(\ell_r,E),q},$$ and
for all $T \in \Pi_{E,p}$
$$\pi_{M(\ell_r,E),q}(T) \le c_p^E \cdot {c_q^{M(\ell_r,E)}}^{-1} \cdot \pi_{E,p}(T).$$
Moreover, if $X$ is a cotype~$2$ space, then for all Banach spaces $Y$
\begin{equation}
\label{inverse}
\Pi_{E,1}(X,Y)=\Pi_{M(\ell_2,E),2}(X,Y).
\end{equation}
\end{lemma}
\proof
 The first inclusion is easy: Let $T:X \rightarrow Y$ be $(E,p)$-summing. Then for
 \mbox{$x_1, \ldots, x_n \in X$} by the
H\"older inequality
\begin{align*}
\norm{(\norm{Tx_k})_1^n}_{M(\ell_r,E)}
&= \sup_{(\lambda_k)_1^n \in B_{\ell_r^n}}\norm{(\lambda_k \cdot\norm{Tx_k})_1^n}_E \\
& \le \pi_{E,p}(T) \cdot c_p^E \cdot \sup_{(\lambda_k)_1^n \in B_{\ell_r^n}}
\sup_{x' \in B_{X'}} \left ( \sum_1^n |\langle x',\lambda_k x_k \rangle|^p \right )^{1/p} \\
&= \pi_{E,p}(T) \cdot c_p^E \cdot   \sup_{x' \in B_{X'}}
\left ( \sum_1^n |\langle x',x_k \rangle|^q \right )^{1/q},
\end{align*}
which gives the claim. The reverse inclusion in the second part follows from the upcoming
Lemma~\ref{multi}:
By a well-known result of Maurey there exists a
constant $C>0$ such that for all $S \in \LL(\ell_\infty^n,X)$:
$\pi_2(S) \le C \cdot \norm{S}$ (see e.\,g. \cite[31.7]{df}). Then
for $T \in \Pi_{M(\ell_2,E),2}(X,Y)$ and $S \in \LL(\ell_\infty^n,X)$ with
$\norm{S}\le 1$ we obtain, together with Lemma~\ref{multi},
$$
\pi_{E,1}(TS) \le \pi_{M(\ell_2,E),2}(T) \cdot \pi_2(S)  \le
C \cdot \pi_{M(\ell_2,E),2}(T),
$$
which by Lemma~\ref{composition} implies $T \in \Pi_{E,1}$.
\qed
\\[15pt]
As announced it remains to prove the following:
\begin{lemma}
\label{multi}
For $1 \le p < q<\infty$ let $1<r<\infty$ such that $1/r=1/p-1/q$. Then
$$
\Pi_{M(\ell_r,E),q} \circ \Pi_r \subset \Pi_{E,p},
$$
and
$$
\pi_{E,p}(TS) \le \pi_{M(\ell_r,E),q}(T) \cdot \pi_r(S)
$$
for $S \in \Pi_r(X,Y)$ and $T \in \Pi_{M(\ell_r,E),q}(Y,Z)$.
\end{lemma}
\proof Let $S $ and $T$ be as in the proposition. By the Pietsch Domination
Theorem there exists a regular Borel probability measure $\mu$ on $B_{X'}$ such that for all $x \in X$
$$
\norm{Sx} \le \pi_r(S) \cdot \left ( \int_{B_{X'}} |\langle x',x \rangle|^r d\mu(x') \right )^{1/r}.
$$
Now take $0 \neq x_1, \ldots,x_n \in X$ and put for $k=1, \ldots,n$
$$
x_k^0 := \left ( \int_{B_{X'}} |\langle x',x_k \rangle|^p d\mu(x') \right )^{-1/r} \cdot x_k.
$$
Then by the H\"older Inequality (and $c_q^{M(\ell_r,E)} \le c_p^E$)
\begin{align*}
\norm{(\norm{TSx_k})_1^n}_E &\le \norm{(\norm{TSx_k^0})_1^n}_{M(\ell_r,E)} \cdot
\left ( \sum_1^n \int_{B_{X'}} |\langle x',x_k \rangle|^p d\mu(x') \right )^{1/r} \\
& \le \pi_{M(\ell_r,E),q}(T) \cdot c_p^E \cdot \sup_{y' \in B_{Y'}}
\left ( \sum_1^n |\langle y',Sx_k^0 \rangle|^q \right )^{1/q} \\
& \qquad \cdot \left ( \sum_1^n \int_{B_{X'}}
 |\langle x',x_k \rangle|^p d\mu(x') \right )^{1/r}.
\end{align*}
Now complete the proof exactly as in \cite{tj70}.
\qed
\par
As in the classical case of $(r,2)$-summing operators, the theory of $(F,2)$-summing operators is
deeply connected to the theory of $s$-numbers. In our case a crucial tool is an extension of an inequality
due to K\"onig, which can be proved exactly as in \cite[2.a.3]{koenig}.
\begin{prop}
\label{weyl}
Let $F$ be a maximal symmetric Banach sequence space such that $\ell_2 \Id F$.
Then $\Pi_{F,2}\subset \ui_F^x$. In particular, for all $T \in \Pi_{F,2}$ and $k$
\begin{equation}
\label{weyleq}
x_k(T) \le \lambda_F(k)^{-1} \cdot c_2^F \cdot \pi_{F,2}(T).
\end{equation}
\end{prop}
The above result allows to give a different proof of the Lemma~\ref{compactm} in the case $p=2$.
\begin{cor}
\label{calkin}
For any Banach space $Y$ any $(F,2)$-summing operator $T: \ell_2 \rightarrow Y$ is compact
whenever $ \ell_2 \Id F \Id c_0$.
\end{cor}
\proof By Proposition~\ref{weyl} we have $\Pi_{F,2}(\ell_2,Y) \subset \ui_F^x(\ell_2,Y) =
\ui_F^a(\ell_2,Y)$ which clearly gives the claim. \qed
\\[15pt]
See Section~6 for the fact that for $2$-convex $F$ the ideals $\Pi_{F,2}$ and
the unitary ideal $\ui_F$ coincide on Hilbert spaces.
 \section{$\pmb{(E,1)}$-summing identity maps}
The well-known results of Bennett \cite{bennett} and Carl \cite{carl}
(proved independently) assure that for $1 \le u \le 2$ the identity map
$\id: \ell_u \Id \ell_2$ is absolutely $(u,1)$-summing. In
 \cite{mastylo} an extension within the setting of Orlicz sequence spaces is presented.
\par
Using interpolation theory we prove as our main result the following proper extension:
\begin{theo}
\label{main}
Let $E$ be a $2$-concave symmetric Banach sequence space. Then the identity map
$\id: E \Id \ell_2$ is $(E,1)$-summing. In other words, for every unconditionally summable
sequence $(x_n)$ in $E$ the scalar sequence $(\norm{x_n}_{\ell_2})$ is contained in $E$.
\end{theo}
The following lemmas are essential:
\begin{lemma}
\label{power}
Let $(E_0,E_1)$ be a relative Calder\'{o}n couple of maximal symmetric Banach sequence spaces
and $E$ an interpolation
space with respect to $(E_0,E_1)$. Then $E^p$  for all $0<p<1$ is an interpolation space with
respect to $(E_0^p,E_1^p)$.
\end{lemma}
\proof It is enough to show that $E^{p}$ is a relative $K$-space with respect
to $(E_{0}^{p}, E_{1}^{p})$, i.e. if whenever $x\in E^{p}$ and
$y\in E_{0}^{p} + E_{1}^{p}$ satisfy
$$
K(t, y; E_{0}^p, E_{1}^p) \leq K(t, x; E_{0}^p, E_{1}^p) \, \, \,
\mbox{for all } t>0,
$$
then it follows that $y\in E^p$.
\par The claim follows from the well-known and easily verified equivalence
for the $K$-functionals, namely
$$
 K(t,x; E_0^{p},E_1^{p}) \asymp
K(t^{1/p},|x|^{1/p};E_0,E_1)^p
$$
for any $x \in E_0^p+E_1^p$ and $t>0$, and the fact that $E$ is a relative $K$-space with
respect to $(E_0, E_1)$.
\qed
\\[15pt]
As an immediate consequence we obtain
\begin{lemma}
\label{interpoldiagonal}
Let $E$ be a maximal symmetric Banach sequence space.
\begin{enumerate}[(a)]
\item If $E$ is $2$-convex, then it is an interpolation space with respect to the couple
$(\ell_2,\ell_\infty)$, i.\,e. there exists an exact interpolation functor $\cal{F}$
such that $E= \mathcal{F}(\ell_2,\ell_\infty)$.
\item
If $E$ is $2$-concave, then $M(\ell_2,E)$ is $2$-convex. In particular, $M(\ell_2,E)$ is an
interpolation space with respect to $(\ell_2,\ell_\infty)$.
\end{enumerate}
\end{lemma}
\proof (a) Without loss of generality we may assume that $\convex{2}(E)=1$. Then $E^2$ is a
 maximal symmetric Banach sequence space, and by Mitiagin~\cite{mitiagin}
(see also \cite[1.b.10]{koenig}) this implies that $E^2$ is an
interpolation space with respect to $(\ell_1,\ell_\infty)$. The claim now follows by the preceding
 lemma and the fact that $(\ell_1,\ell_\infty)$ is a relative Calder\'{o}n couple
 (see e.\,g. \cite[2.6.9]{bk}). \\
 (b) Without loss of generality we may assume that $\concave{2}(E)=1$. Then $E^\times$ is
 $2$-convex with $\convex{2}(E)=1$, hence $(E^\times)^2$ and therefore also
 $((E^\times)^2)^\times$ are normed.
Consequently, $M(\ell_2,E)=(((E^\times)^2)^\times)^{1/2}$ is $2$-convex.
\qed
\\[15pt]
For the sake of
 completeness  we give a proof of the following easy and well-known result:
\begin{lemma}
\label{vector}
Let $E$ and $F$ be Banach sequence spaces, $X$ a Banach space and
${\cal{F}}$ an exact interpolation functor. Then
$$
\norm{\id : {\cal{F}}(E_n(X),F_n(X)) \Id {\cal{F}}(E_n,F_n)(X)} \le 1.
$$
\end{lemma}
\proof For any given $x_1, \ldots,x_n \in X$ let $x_1', \ldots, x_n' \in X'$ be such that
$\norm{x_i'}=1$ and \mbox{$\langle x_i',x_i \rangle=\norm{x_i}$.} Then for
$T: \R^n(X) \rightarrow \R^n$
defined by $T((y_i)_1^n):= (\langle x_i',y_i \rangle)_1^n$ we obviously have $\norm{T: E_n(X) \rightarrow
E_n} \le 1$ and $\norm{T: F_n(X) \rightarrow
F_n} \le 1$, hence
$$
\norm{T: {\cal{F}}(E_n(X),F_n(X)) \rightarrow {\cal{F}}(E_n,F_n)} \le 1.
$$
\begin{samepage}
Thus
\begin{align*}
\norm{(x_i)_1^n}_{{\cal{F}}(E_n,F_n)(X)} &=
\norm{(\norm{x_i})_1^n}_{{\cal{F}}(E_n,F_n)} \\
& = \norm{(\langle x_i',x_i \rangle)_1^n}_{{\cal{F}}(E_n,F_n)} \\
& \le \norm{(x_i)_1^n}_{{\cal{F}}(E_n(X),F_n(X))}.
\end{align*}
\qed
\end{samepage}
\par The following lemma partially extends (in the lattice case) results of Pisier and Kouba on
the complex interpolation
 of spaces of operators (see \cite{kouba}, \cite{pisier90} and also \cite{dm99}).
 Recall that $\ell_2 = M(\ell_2,\ell_1)$ and $\ell_\infty=M(\ell_2,\ell_2)$;
 then the statement below says that under the given assumption the interpolation property of the
 spaces of multipliers (diagonal operators) can be transferred into the corresponding interpolation property
 of the associated spaces of bounded operators (at least in the finite-dimensional case). Note
 that a formula for the reverse inclusion holds whenever ${\cal{F}}(\ell_1,\ell_2) \Id E$.
\begin{lemma}
\label{kouba}
For a $2$-concave symmetric Banach sequence space $E$ let ${\cal{F}}$ be an exact interpolation
functor such that $M(\ell_2,E) \Id {\cal{F}}(\ell_2,\ell_\infty)$. Then
\begin{equation}
\label{koubaformula}
\sup_{m,n} \norm{\id: \LL(\ell_2^m,E_n) \Id {\cal{F}}(\LL(\ell_2^m,\ell_1^n),
\LL(\ell_2^m,\ell_2^n))}
\le \sqrt{2} \cdot c_{M(\ell_2,E)}^{\mathcal{F}(\ell_2,\ell_\infty)} \cdot \concave{2}(E).
\end{equation}
\end{lemma}
\proof Let $T \in \LL(\ell_2^m,E_n)$. By a variant of the Maurey--Rosenthal Factorization Theorem
(see
\cite[4.2]{defant} and also \cite{lust}) there exist an operator $R \in \LL(\ell_2^m,\ell_2^n)$ and $\lambda
\in \R^n$ such that
$$
\norm{R} \cdot \norm{\lambda}_{M(\ell_2^n,E_n)} \le \sqrt{2} \cdot \concave{2}(E) \cdot
\norm{T}
$$
and $T$ factorizes as follows:
\begin{center}
\resetparms
\Vtriangle[\ell_2^m`E_n`\ell_2^n;
T`R`M_\lambda].
\end{center}
Obviously the map $\Phi$ defined by
$$
\Phi (\mu) := M_{\mu} \circ R, \quad \mu \in \R^n
$$
 maps the couple
$(\ell_{2}^{n}, \ell_{\infty}^{n})$ into the couple
$
({\cal L}(\ell_{2}^{m}, \ell_{1}^{n}),
{\cal L}(\ell_{2}^{m}, \ell_{2}^{n}))
$
such that both restrictions have norm less or equal $\|R\|$. Hence by the interpolation property
and
the assumption
 the restriction map
$$
\Phi: M(\ell_2^n,E_n) \rightarrow {\cal{F}}(\LL(\ell_2^m,\ell_1^n),\LL(\ell_2^m,\ell_2^n))
$$
has norm $\le c_{M(\ell_2,E)}^{\mathcal{F}(\ell_2,\ell_\infty)} \cdot\norm{R}$. Thus we obtain
\begin{samepage}
\begin{align*}
\norm{T}_{{\cal{F}}(\LL(\ell_2^m,\ell_1^n),\LL(\ell_2^m,\ell_2^n))} & =
\norm{M_\lambda \circ R}_{{\cal{F}}(\LL(\ell_2^m,\ell_1^n),\LL(\ell_2^m,\ell_2^n))} \\
& \le c_{M(\ell_2,E)}^{\mathcal{F}(\ell_2,\ell_\infty)} \cdot \norm{R} \cdot \norm{\lambda}_{M(\ell_2^n,E_n)} \\
& \le \sqrt{2} \cdot c_{M(\ell_2,E)}^{\mathcal{F}(\ell_2,\ell_\infty)} \cdot \concave{2}(E) \cdot \norm{T}_{\LL(\ell_2^m,E_n)}.
\end{align*}
\qed
\end{samepage}
\\[20pt]
Now we are ready to give a proof of Theorem~\ref{main}: According to Lemma~\ref{interpoldiagonal}
let ${\cal{F}}$ be an interpolation
functor with $M(\ell_2,E)={\cal{F}}(\ell_2,\ell_\infty)$.
We consider the mapping
$$
\Phi^{m,n}: ({\cal L}(\ell_{2}^{m},\ell_{1}^{n}),
{\cal L}(\ell_{2}^{m}, \ell_{2}^{n})) \to
(\ell_{2}^{m}(\ell_{2}^{n}), \ell_{\infty}^{m}(\ell_{2}^{n}))
$$
defined by $\Phi^{m,n}(S):= (Se_i)_{1}^{m}$. By \eqref{discrete} we have
$$
\sup_{m} \|\Phi^{m,n}: {\cal L}(\ell_{2}^{m}, \ell_{1}^{n})
\to \ell_{2}^{m}(\ell_{2}^{n})\| = \pi_{2}(\mbox{id}: \ell_{1}^{n}
\hookrightarrow \ell_{2}^{n}\| =1
$$
and
$$
\sup_{m} \|\Phi^{m,n}: {\cal L}(\ell_{2}^{m}, \ell_{2}^{n})
\to \ell_{\infty}^{m}(\ell_{2}^{n})\| = \|\mbox{id}: \ell_{2}^{n}
\hookrightarrow \ell_{2}^{n}\| = 1.
$$
Then by the interpolation property we obtain that
$$
\Phi^{m,n}: {\cal{F}}(\LL(\ell_2^m,\ell_1^n),\LL(\ell_2^m,\ell_2^n)) \rightarrow
{\cal{F}}(\ell_2^m(\ell_2^n),\ell_\infty^m(\ell_2^n))
$$
also has norm $\le 1$. Now by the preceding lemma
$$
\norm{\Phi^{m,n}: \LL(\ell_2^m,E_n) \rightarrow M(\ell_2^m, E_m)(\ell_2^n)} \le \sqrt{2} \cdot
c_{M(\ell_2,E)}^{\mathcal{F}(\ell_2,\ell_\infty)} \cdot  \concave{2}(E).
$$
Hence, since $\sup_m \norm{\Phi^{m,n}} = \pi_{M(\ell_2,E),2}(\id: E_n \Id \ell_2^n)$,
$$
\pi_{M(\ell_2,E),2}(\id : E_n \Id \ell_2^n) \le \sqrt{2} \cdot
c_{M(\ell_2,E)}^{\mathcal{F}(\ell_2,\ell_\infty)}
\cdot \concave{2}(E),
$$
and since $\bigcup_n E_n$ is dense in $E$, this implies $(\id : E \Id \ell_2) \in
\Pi_{M(\ell_2,E),2}$. The final statement then follows from \eqref{inverse}. \qed
\par
Theorem~\ref{main} is best possible in the following sense:
\begin{cor}
\label{equivalent}
Let $E$ and $F$ be $2$-concave  symmetric Banach sequence spaces. Then
\begin{equation}
\pi_{F,1}(\id: E_n \Id \ell_2^n) \asymp \norm{\id: E_n \Id F_n}.
\end{equation}
In particular, $\id: E \Id \ell_2$ is $(F,1)$-summing if and only if $E \Id F$.
\end{cor}
\proof The upper estimate follows from Theorem~\ref{main} by factorization; for the lower estimate
 we may assume without loss of generality that $\concave{2}(E)=\concave{2}(F)=1$. Observe that
 for $\lambda \in \R^n$ one has $\norm{\lambda}_{M(\ell_2^n,F_n)} \le
 \pi_{M(\ell_2,F),2}(M_\lambda: \ell_2^n \rightarrow \ell_2^n)$ (simply take in the definition
 of $(M(\ell_2,F),2)$-summing $x_i=e_i$), hence, by Lemma~\ref{composition} and
 Lemma~\ref{inclusion} as well as
 \eqref{diagonalpower},
 \begin{samepage}
\begin{align*}
\pi_{F,1}(\id:E_n \Id \ell_2^n) & \ge \pi_{M(\ell_2,F),2}(\id:E_n \Id \ell_2^n) \\
&= \sup_m \sup_{\norm{S: \ell_2^m \rightarrow E_n} \le 1}
\pi_{M(\ell_2,F),2}(\ell_2^m \xrightarrow{S} E_n \xrightarrow{\id} \ell_2^n) \\
& \ge \sup_{\norm{\lambda}_{M(\ell_2^n,E_n)} \le 1} \pi_{M(\ell_2,F),2}(M_\lambda: \ell_2^n
\rightarrow \ell_2^n) \\
& \ge \sup_{\norm{\lambda}_{M(\ell_2^n,E_n)} \le 1} \norm{\lambda}_{M(\ell_2^n,F_n)} \\
& = \norm{\id : M(\ell_2^n, E_n) \Id M(\ell_2^n,F_n)} \\
& = \norm{\id : (((E_n^\times)^2)^\times)^{1/2} \Id (((F_n^\times)^2)^\times)^{1/2}} \\
& = \norm{\id : E_n \Id F_n}.
\end{align*}
\qed
\end{samepage}
\par
As a counterpart to Corollary~\ref{equivalent} we show that in Theorem~\ref{main} the Hilbert
space $\ell_2$ is minimal in the following sense:
\begin{cor}
\label{bestpossible}
Let $E$ and $F$ be maximal symmetric Banach sequence where E is $2$-concave. Then
\begin{equation}
\pi_{E,1}(\id: E_n \Id F_n) \asymp \norm{\ell_2^n \Id F_n}.
\end{equation}
In particular, $\id: E \Id F$ is $(E,1)$-summing if and only if $\ell_2 \Id F$.
\end{cor}
\proof Again the upper estimate obviously follows by factorization from Theorem~\ref{main}. For the
lower estimate note that by \cite[p.~237 (3)]{CD97} (which is also valid for
$\lfloor n/2 \rfloor +1$ instead of $[n/2]$)
$$
x_{\lfloor n/2 \rfloor +1}(\id : E_n \Id F_n) \ge \frac{1}{\sqrt{2}} \cdot
\frac{\norm{\id :\ell_2^n \Id F_n}}{\norm{\id: \ell_2^n \Id E_n}},
$$
hence, by Lemma~\ref{inclusion}, \eqref{weyleq} and \cite[p.~237 (1)]{CD97},
\begin{samepage}
\begin{align*}
\pi_{E,1}(\id : E_n \Id F_n) & \ge \pi_{M(\ell_2,E),2}(\id : E_n \Id F_n) \\
& \ge \frac{ \norm{\id: \ell_2^n \Id F_n}}{\sqrt{2}} \cdot \frac{\norm{\id : \ell_2^{\lfloor n/2 \rfloor+1} \Id
E_{\lfloor n/2 \rfloor+1}}}{\norm{\id: \ell_2^n \Id E_n}}  \\
&\ge \frac{\norm{\id: \ell_2^n \Id F_n}}{\sqrt{2}} \cdot \frac{\norm{\id : \ell_2^{n-\lfloor n/2 \rfloor} \Id
E_{n-\lfloor n/2 \rfloor}}}{\norm{\id: \ell_2^n \Id E_n}}   \\
& \ge \frac{\norm{\id: \ell_2^n \Id F_n}}{\sqrt{2}} \cdot
\frac{a_{\lfloor n/2 \rfloor +1}(\id: \ell_2^n \Id E_n)}{\norm{\id: \ell_2^n \Id E_n}}
 \\
&\ge \frac{\norm{\id: \ell_2^n \Id F_n}}{2}.
\end{align*}
\qed
\end{samepage}
\par We note that in general the assumption that a symmetric sequence
space $E$ is 2-concave is essential in Theorem~\ref{main}, even in the class
of Orlicz sequence spaces. This follows from the following proposition
and the fact that there is an example, constructed by Kalton \cite{kalton}
(see also \cite[4.c.3]{lt77}), of an Orlicz sequence space $\ell_{\varphi}$
such that the identity map $\mbox{id}: \ell_{\varphi} \hookrightarrow \ell_{2}$
is not a strictly singular operator, i.e. $\mbox{id}$ is an isomorphism
on some infinite dimensional closed subspace of $\ell_{\varphi}$.
\begin{prop}
\label{ss1}
Let $E \Id \ell_2$ be a Banach sequence
space not equivalent to $\ell_2$. Then  the identity map $\id: E \Id \ell_2$ is strictly singular
 whenever it is $(E,1)$-summing.
\end{prop}
\proof
Suppose that $\mbox{id}: E \hookrightarrow \ell_2$ is not strictly singular.
Thus there exists an infinite dimensional closed subspace
$X$ of $E$ such that the restriction of $\mbox{id}$ to $X$
is an isomorphism from $X$ into $\ell_{2}$. Let $P: \ell_{2} \to X$ be
a continuous linear projection. By assumption
$\mbox{id}: E \hookrightarrow \ell_{2}$ is $(E, 1)$-summing,
and thus by Lemma~\ref{inclusion}
$T=\mbox{id}\circ P: \ell_2 \to \ell_2$ is $(M(\ell_2, E), 2)$-summing.
Since $E \neq \ell_{2}$, we get $M(\ell_{2}, E) \hookrightarrow c_0$
(see \eqref{czero}). An application of Lemma~\ref{compact} yields that $T$ is compact which
  contradicts the fact that $T$ on $X$ is the identity.
\qed
\\[15pt] In view of Theorem~\ref{main}
 the following trivial consequence  seems to be of independent
interest.
\begin{cor}
\label{ss2}
If $E$ is a symmetric Banach sequence space not
equivalent to $\ell_2$ and such that the inclusion map
$\mbox{id}: E \hookrightarrow \ell_2$  is not strictly singular,
then $E$ does not have cotype $2$.
\\[15pt]
\em Combining Proposition~\ref{ss1} and Corollary~\ref{ss2} we see that Kalton's example $\ell_\varphi$
is not $2$-concave and $\id: \ell_\varphi \Id \ell_2$ is not $(\ell_\varphi,1)$-summing.
\end{cor}

\section{Applications to approximation numbers of identity operators}
Of
special interest for applications (e.\,g. in approximation theory)
are formulas for the asymptotic behavior of approximation numbers of
finite-dimensional identity
operators. One of the first well-known results in this direction
is due to Pietsch \cite{pietsch74}: For $1 \le k \le n$ and $1 \le p<q \le \infty$
\begin{equation}
\label{lpformula} a_k(\id: \ell_q^n \Id \ell_p^n) = (n-k+1)^{1/p-1/q}.
\end{equation}
For the special case $1 \le p < q=2$ let us rewrite this as follows:
\begin{equation}
\label{l2formula} a_k(\id: \ell_2^n \Id \ell_p^n) = \frac{\lambda_{\ell_p}(n-k+1)}{(n-k+1)^{1/2}}.
\end{equation}
 Using Theorem~\ref{main}
we show this formula---at
least asymptotically---for all $2$-concave symmetric Banach sequence spaces $E$
instead of $\ell_p$:
\begin{theo}
\label{formula1}
Let $E$ be a $2$-concave symmetric Banach sequence space. Then for all \mbox{$1 \le k \le n$}
\begin{equation}
\label{formula2}
a_k(\id: \ell_2^n \Id E_n) \asymp \frac{\lambda_E(n-k+1)}{(n-k+1)^{1/2}}.
\end{equation}
\end{theo}
\noindent
The proof needs the special case $k=1$, a result due to Szarek and
Tomczak-Jaegermann~\cite[Proposition~2.2]{stj}:
Under the assumption of the theorem
\begin{equation}
\label{stj}
a_1(\id : \ell_2^n \Id E_n) = \norm{\id: \ell_2^n \Id E_n} \asymp
\frac{\lambda_E(n)}{n^{1/2}}.
\end{equation}
Proof of Theorem~\ref{formula1}:
First we claim that it is enough to show
\begin{equation}
\label{formula3}
a_k(\id: \ell_2^n \Id E_n) \asymp \norm{\textstyle{\sum_1^{n-k+1} e_i}}_{M(\ell_2,E)}.
\end{equation}
Indeed, the right hand side in \eqref{formula3} is obviously equal to $\norm{\id: \ell_2^{n-k+1}
 \Id E_{n-k+1}}$, and by \eqref{stj} this is
 asymptotically equivalent to the right hand side in \eqref{formula2}.
\\[10pt]
The upper estimate in \eqref{formula3} is straightforward:
Put $\lambda:= \sum_1^{n-k+1} e_i \in \K^n$ and \mbox{$\mu:= \sum_{n-k+2}^n e_i \in \K^n$.}
Since the diagonal operator $M_\mu: \ell_2^n \rightarrow E_n$ has rank $k-1$, we obtain
$$
a_k(\id: \ell_2^n \Id E_n) \le \norm{\id-M_\mu}=\norm{M_\lambda}=
\norm{\textstyle{\sum_1^{n-k+1} e_i}}_{M(\ell_2,E)}.
$$
On the other hand, by a result of \cite{CD92}
$$
a_k(\id:\ell_2^n \Id E_n)=x_{n-k+1}(\id:E_n \Id \ell_2^n)^{-1},
$$
so that the lower estimate in \eqref{formula3} follows from
\begin{equation}
\label{formula4}
x_k(\id:E_n \Id \ell_2^n) \prec \norm{\textstyle{\sum_1^k e_i}}_{M(\ell_2,E)}^{-1}.
\end{equation}
In order to check \eqref{formula4} note that by Theorem~\ref{main} the identity map
$\id: E \Id \ell_2$ is $(E,1)$-summing. Hence by
Lemma~\ref{inclusion} it is also $(M(\ell_2,E),2)$-summing, and by the generalized
K\"onig inequality~\eqref{weyleq} we obtain
\begin{align*}
x_k(\id : E_n \Id \ell_2^n) &\le \norm{\textstyle{\sum_1^k e_i}}_{M(\ell_2,E)}^{-1} \cdot
\pi_{M(\ell_2,E),2}(\id: E_n \Id \ell_2^n) \\
&\le  \norm{\textstyle{\sum_1^k e_i}}_{M(\ell_2,E)}^{-1} \cdot
\pi_{M(\ell_2,E),2}(\id: E \Id \ell_2),
\end{align*}
which completes the proof. \qed \\[10pt]
To illustrate formula~\eqref{formula2} we consider Lorentz and Orlicz sequence spaces.
\begin{samepage}
\begin{cor}
\begin{enumerate}[(a)]
\item
Let $1 <p<2$ and $1 \le q \le 2$. Then for all $1 \le k \le n$
\begin{equation}
\label{lorentz}
a_k( \id : \ell_2^n \Id \ell_{p,q}^n) \asymp (n-k+1)^{1/p-1/2}.
\end{equation}
\item
\label{lorentz2}
Let $1<p<2$ and $w$ be a Lorentz sequence such that $n \cdot w_n^{2/(2-p)} \asymp \sum_1^n
w_i^{2/(2-p)}$. Then for all $1 \le k \le n$
\begin{equation}
a_k(\id: \ell_2^n \Id d_n(w,p)) \asymp (n-k+1)^{1/p-1/2} \cdot w_{n-k+1}^{1/p}.
\end{equation}
\item
\label{orlicz}
Let $\varphi$ be an Orlicz function such that the function $t \mapsto \varphi(\sqrt{t})$ is
equivalent to a concave function. Then
\begin{equation}
a_k(\id : \ell_2^n \Id \ell_\varphi^n) \asymp \frac{(\varphi^{-1}(1/(n-k+1)))^{-1}}{(n-k+1)^{1/2}}.
\end{equation}
 \end{enumerate}
\em Note that (a) is---asymptotical---the same result as for $\ell_p$ (see \eqref{lpformula});
although $\ell_{p,q}$ is
 ''very close`` to $\ell_p$, one may have expected an additional logarithmic term.
\end{cor}
\end{samepage}
\proof By Theorem~\ref{formula1} it is enough to ensure that all spaces considered in the corollary are
$2$-concave. For the Lorentz sequence spaces $\ell_{p,q}$ this is due to Creekmore~\cite{creekmore}
(see e.\,g. also \cite{defant}), for
the Lorentz sequence spaces $d(w,p)$ see Reisner~\cite{reisner}, and for Orlicz sequence spaces this
is contained in \cite{komarchev}.\qed
\section{Applications to eigenvalues of compact operators and unitary ideals}
By Pitt's Theorem every operator $T$ on $\ell_2$ with values in $\ell_u$, $1 \le u <2$, is
compact. The original Bennett--Carl result implies (see e.\,g. \cite[2.b.11]{koenig}) that its
sequence of singular numbers is contained in $\ell_r$, $1/r=1/u-1/2$, and by Weyl's Inequality
(see e.\,g. \cite[1.b.9]{koenig}) even its sequence of eigenvalues is contained in $\ell_r$.
Weyl's inequality also holds for arbitrary maximal symmetric Banach sequence spaces: If the
singular numbers of a compact operator on a Hilbert space is contained in a certain maximal symmetric
sequence space $F$, then the same is true for its sequence of eigenvalues (see \cite[1.b.10]{koenig}).
 Together with Theorem~\ref{main} this implies the following extension of the result mentioned
 above:
\begin{theo}
\label{compact}
Let $E \hookrightarrow \ell_2$ be a $2$-concave symmetric Banach sequence space not equivalent
 to $\ell_2$,
and $T \in \LL(\ell_2,\ell_2)$
be an operator with values in $E$. Then $T \in \ui_{M(\ell_2,E)}$. In particular, its
sequence of eigenvalues   $(\lambda_n(T))$ is contained in $M(\ell_2,E)$.
\end{theo}
\proof The assumption $E \not= \ell_2$ together with Corollary~\ref{equivalent} assures that
 the identity operator on $\ell_2$ is not contained in $\Pi_{M(\ell_2,E),2}$, and by a result
 of Calkin (see \cite[2.11.11]{eigenvalues}) it follows that every operator in
 $\Pi_{M(\ell_2,E),2}(\ell_2,\ell_2)$ is compact; alternatively one may directly use
 Lemma~\ref{compactm} together with \eqref{czero}. Now by Theorem~\ref{main} and the ideal property
 the operator $T : \ell_2 \xrightarrow{T} E \xrightarrow{\id} \ell_2$ is contained in
 $\Pi_{M(\ell_2,E),2}(\ell_2,\ell_2)$ and therefore compact, and by Proposition~\ref{weyl}
 the sequence of Weyl (=singular) numbers $(x_n(T))$ is contained in $M(\ell_2,E)$. The second
 claim now follows by Weyl's inequality mentioned above. \qed

\par
Next we discuss an alternative approach to Theorem~\ref{main} using interpolation
of unitary ideals; we first
illustrate our idea by considering the original result of Bennett and Carl:
\\[10pt]
Let $1\le u <2$. By Lemma~\ref{composition} and \eqref{inverse} the identity map
$I_u: \ell_u \Id \ell_2$ is absolutely $(u,1)$-summing whenever the composition $I_uS$
for any operator
$S: \ell_2 \rightarrow \ell_u$  is absolutely $(r,2)$-summing $(1/r=1/u-1/2)$.
By the (classical) Maurey--Rosenthal
Factorization Theorem there exist an operator $R \in \LL(\ell_2,\ell_2)$ and
 $\lambda \in \ell_r$ such that
 $S$ factorizes as follows:
\begin{center}
\resetparms
\Vtriangle[\ell_2`\ell_u`\ell_2;
S`R`M_\lambda].
\end{center}
Then obviously the operator $I_u M_\lambda: \ell_2 \rightarrow \ell_2$ is contained in the
Schatten-$r$-class $\ui_r$. By a
result of Mitiagin (see e.\,g. \cite[10.3]{djt}) $\ui_r = \Pi_{r,2}(\ell_2,\ell_2)$,
hence $I_u S = I_u M_\lambda R$ is absolutely $(r,2)$-summing which gives the claim. \qed
\\[10pt]
Mitiagin's result and its proof are of interpolative nature.
Alternatively, the inclusion \mbox{$\ui_r \subset \Pi_{r,2}(\ell_2,\ell_2)$} can be proved by complex interpolation of the
border cases $\Pi_{2,2}(\ell_2,\ell_2)=\ui_2$ and
$\Pi_{\infty,2}(\ell_2,\ell_2)=\LL(\ell_2,\ell_2)=\ui_\infty$: For $\theta:=2/r$
$$
\ui_r =[\ui_2,\ui_\infty]_\theta =[\Pi_{2,2}(\ell_2,\ell_2),\Pi_{\infty,2}(\ell_2,\ell_2)]_\theta
 \subset \Pi_{r,2}(\ell_2,\ell_2).
$$
The starting point for our alternative approach to Theorem~\ref{main} now is an
extension of Mitiagin's result for which we  need
 the following generalization of a result due to K\"onig (cf. \cite[2.c.10]{koenig}).
 \begin{samepage}
\begin{lemma}
\label{lemma35}
Let ${\cal F}$ be an interpolation functor
and $(E_0, E_1)$  a couple of Banach sequence spaces with
$\ell_{p} \hookrightarrow E_j$ and $\norm{e_n}_{E_j}=1$ for all $n$, $j=0,1$. Then for arbitrary
Banach spaces $X$ and $Y$, we have
$$
{\cal F} \big (\Pi_{E_0, p}(X, Y), \Pi_{E_1, p}(X,Y)\big)
\hookrightarrow \Pi_{{\cal F}(E_0, E_1), p}(X, Y).
$$
\end{lemma}
\end{samepage}
\proof For fixed vectors $x_1, \ldots, x_n \in X$ with $\sup_{x' \in B_{X'}}\sum_{j=1}^n
|\langle x',x_j \rangle|^p  \le 1$, we
define \mbox{$\Phi (T):= (Tx_j)_{j=1}^{n}$} for $T\in {\cal L}(X, Y)$.
Clearly
$$
\Phi : \big (\Pi_{E_0,p}(X, Y), \Pi_{E_1,p}(X,Y)\big)
\to (E_{0n}(Y), E_{1n}(Y))
$$
with norm $\le 1$.
By interpolation and Lemma~\ref{vector} we obtain that
$$
\Phi: {\cal F} \big (\Pi_{E_0, p}(X, Y), \Pi_{E_1, p}(X,Y)\big)
\to {\cal F}(E_{0n}, E_{1n})(Y)
$$
with norm $\le 1$. This yields that
$$
\pi_{{\cal F}(E_0, E_1),p}(T) \leq \norm{T}_{{\cal F}(\Pi_{E_0,p}(X,Y), \Pi_{E_1,p}(X,Y))}
$$
for any $T\in {\cal F}(\Pi_{E_0,p}(X,Y), \Pi_{E_1,p}(X,Y))$.
\qed
\\[10pt] The following theorem now admits the announced alternative proof of Theorem~\ref{main}
 exactly
 as it was done above for the original Bennett--Carl result---but it also seems to be of
 independent interest.
\begin{theo}
\label{interpolunitary}
Let $F$ be a $2$-convex maximal symmetric Banach sequence space. Then $\Pi_{F,2}(\ell_2,\ell_2)=
\ui_F$.
\end{theo}
\proof The inclusion $\Pi_{F,2}(\ell_{2}, \ell_2) \hookrightarrow
{\cal S}_{F}$ is contained in Proposition \ref{weyl}. For the reverse inclusion,
we note that if $E \hookrightarrow c_{0}$ is a maximal symmetric space, then
$E$ is an interpolation space with respect to $(\ell_{1}, c_{0})$.
Assume without loss of generality that $F \neq \ell_{\infty}$ and
\mbox{${\bf M_{(2)}}(F)=1$.} By the symmetry of $F$, it follows that
$F\hookrightarrow c_{0}$. In consequence $F^{2}$ is an interpolation
space with respect to $(\ell_{1}, c_{0})$, and thus by Lemma \ref{power}, $F$ is
an interpolation space with respect to $(\ell_{2}, c_0)$ (note that $(\ell_1,c_0)$ is a
relative Calder\'{o}n couple since $(\ell_1,\ell_\infty)$ is). Hence
there exists an exact interpolation functor ${\cal F}$ such that
$F={\cal F}(\ell_{2}, c_0)$. By applying Lemma \ref{lemma35}, we obtain
$$
{\cal F}(\Pi_{\ell_2, 2}(\ell_2, \ell_2),
\Pi_{c_0, 2}(\ell_{2}, \ell_{2})) \hookrightarrow \Pi_{F,2}(\ell_2, \ell_2)
$$
The claim now follows by the fact that ${\cal K}(\ell_{2}) \hookrightarrow
\Pi_{c_0, 2}(\ell_2, \ell_2)$ (${\cal K}(\ell_{2})$ denotes the space of compact
operators on $\ell_{2}$) and by a result on interpolation of unitary
ideals due to Arazy [Ara78]:
${\cal S}_{F} = {\cal S}_{{\cal F}(\ell_2, c_0)} =
{\cal F}({\cal S}_2, {\cal K}(\ell_2))$.
 \qed
\par Another nice application of Theorem~\ref{interpolunitary} is the following:
 \begin{cor}
Let $F$ be a maximal symmetric Banach sequence space such that
\mbox{$\ell_2 \Id F$.} Then
for every Banach space $X$ with $\dim X=n$
\begin{equation}
\label{id}
 \pi_{F,2}(\id_X) \ge C^{-1} \cdot \lambda_F(n) ,
\end{equation}
where $C>0$ is a constant depending on $F$ only. Moreover, if $F$ is $2$-convex, then even
\begin{equation}
\label{id2}
C^{-1} \cdot \lambda_F(n) \le \pi_{F,2}(\id_X) \le C \cdot \lambda_F(n).
\end{equation}
\end{cor}
\proof Let $G:=m_{\lambda_F}$ be the Marcinkiewicz sequence space associated to $F$.
By Proposition~\ref{weyl} and the fact that the continuous
inclusion $F \Id G$ is of norm one we get that
$$
\norm{(x_k(\id_X))_1^n}_G \le c_2^F \cdot \pi_{F,2}(\id_X).
$$
Since $G$ is a maximal symmetric Banach sequence space, it follows by
 the generalized Weyl Inequality \cite[2.a.8]{koenig} that
\begin{align*}
\norm{\textstyle{\sum_1^n e_i}}_F &= \sup_{1 \le k\le n} \norm{\textstyle{\sum_1^k e_i}}_F
\cdot 1 =
\norm{(\lambda_k(\id_X))_1^n}_G \\
&\le 2\sqrt{2e} \cdot \norm{(x_k(\id_X))_1^n}_G \le 2 \sqrt{2e} \cdot c_2^F \cdot \pi_{F,2}(\id_X),
\end{align*}
where $\lambda_k(\id_X)$ is the $k$-th eigenvalue of $\id_X$.
For the reverse estimate note that for an operator $T: Y \rightarrow Z$ of rank $n$ one has
$$
\pi_{F,2}(T)=\sup \{\pi_{F,2}(TS) \, | \,S \in \LL(\ell_2^n,Y), \quad \norm{S} \le 1 \}
$$
(check the proof of \cite[11.3 and 9.7]{tj}). Now let $S \in \LL(\ell_2^n,X)$. Then by
Theorem~\ref{interpolunitary}
$$
\pi_{F,2}(\id_X \circ S) \le \pi_{F,2}(\id_{\ell_2^n}) \cdot \norm{S} \le \tilde{C} \cdot
\norm{\textstyle{\sum_1^n e_i}}_F
\cdot \norm{S},
$$
where $\tilde{C}>0$ is a constant only depending on $F$.
\qed
\section{Complex interpolation in the range}
Based on the case $v=2$, Bennett and Carl also proved that for $1 \le u \le v \le 2$ the identity operator $\id: \ell_u
\Id \ell_v$ is absolutely $(r,2)$-summing whenever $1/r=1/u-1/v$. By using Theorem~\ref{main} and
complex interpolation in the range we obtain the following formal extension of our main result:
\begin{prop}
\label{complex}
Let $E$ be a $2$-concave symmetric Banach sequence space. Then for \mbox{$0 \le \theta<1$} the
identity operator $\id : E \hookrightarrow [\ell_2,E]_\theta$ is absolutely
$(M(\ell_2,E)^{1-\theta},2)$-summing.
\end{prop}
This now enables us to give an extension of the original Bennett--Carl result within the
framework of Lorentz sequence spaces:
\begin{cor}
\label{lorentzbc}
Let $1 <u_1 <v_1<2$ and $1\le u_2 \le v_2 \le 2$ be such that either $u_2=v_2=2$ or
$\frac{1/{v_1}-1/2}{1/{u_1}-1/2} = \frac{1/{v_2}-1/2}{1/{u_2}-1/2}$. Then
$$( \id: \ell_{u_1,u_2} \Id \ell_{v_1,v_2}) \in \Pi_{\ell_{r_1,r_2},2},$$
where $1/{r_1}=1/{u_1}-1/{v_1}$ and $1/{r_2}=1/{u_2}-1/{v_2}$.
\end{cor}
\proof This directly follows from the preceding proposition and the fact that for
$\theta:=\frac{1/{v_1}-1/2}{1/{u_1}-1/2}$ by the reiteration
theorem \cite[4.7.2]{BL} one has
$
[\ell_2, \ell_{u_1,u_2}]_\theta = \ell_{v_1,v_2};
$
finally,
$$M(\ell_2,\ell_{u_1,u_2})^{1-\theta}=\ell_{\tilde{u}_1,\tilde{u}_2}^{1-\theta}=\ell_{r_1,r_2},$$
where $1/{\tilde{u}_1}=1/{u_1}-1/2$ and $1/{\tilde{u}_2}=1/{u_2}-1/2$. \qed
\par
For our applications of this result we need the following two statements:
\begin{lemma}
Let $F$ be a maximal symmetric sequence space such that $\ell_2 \Id F$. Then for every invertible
operator $T: X \rightarrow Y$ between two $n$-dimensional Banach spaces and all $1 \le k \le n$
\begin{equation}
\label{gelfand}
c_k(T) \ge C^{-1} \cdot \frac{\lambda_F(n-k+1)}{\pi_{F,2}(T^{-1})},
\end{equation}
where $C:= 2 \sqrt{2e} \cdot c_2^F$.
\end{lemma}
\proof We copy the proof of \cite[p.~231]{CD97} for the $2$-summing norm. Take a subspace
$M \subset X$ with $\codim \, M <k$. Then
$$
n-k+1 \le \dim M,
$$
hence by \eqref{id}
$$
\norm{\textstyle{\sum_1^{n-k+1} e_i}}_F \le \norm{\textstyle{\sum_1^{\dim M} e_i}}_F \le
C \cdot \pi_{F,2}(\id_M).
$$
Clearly (by the injectivity of $\Pi_{F,2}$)
$$
\pi_{F,2}(\id_M) = \pi_{F,2}(\id: M \hookrightarrow X),
$$
therefore the commutative diagram
\begin{center}
\resetparms
\Vtriangle[M`X`Y;\id`T_{|M}`T^{-1}]
\end{center}
gives, as desired,
$
\norm{\textstyle{\sum_1^{n-k+1} e_i}}_F \le \norm{T_{|M}} \cdot C \cdot \pi_{F,2}(T^{-1}).
$
\qed
\\[15pt]
The following two results extend \eqref{formula2} and \eqref{lorentz}:
\begin{prop}
Let $E$ be a $2$-concave symmetric Banach sequence space. Then for \mbox{$0 \le \theta<1$} and all
$1 \le k \le n$
\begin{equation}
\label{ckinterpol}
a_k(\id: [\ell_2^n,E_n]_\theta \Id E_n) \asymp c_k(\id: [\ell_2^n,E_n]_\theta \Id E_n) \asymp
\left (\frac{\lambda_E(n-k+1)}{(n-k+1)^{1/2}}\right )^{1-\theta}.
\end{equation}
\end{prop}
\proof The estimate
$$
c_k(\id: [\ell_2^n,E_n]_\theta \Id E_n) \succ \left (
\frac{\lambda_E(n-k+1)}{(n-k+1)^{1/2}} \right )^{1-\theta}
$$
follows from Proposition~\ref{complex} and \eqref{gelfand} together with \eqref{stj}. Obviously
$$
c_k(\id: [\ell_2^n,E_n]_\theta \Id E_n) \le a_k(\id: [\ell_2^n,E_n]_\theta \Id E_n) \le
\norm{\id : [\ell_2^{n-k+1},E_{n-k+1}]_\theta \Id E_{n-k+1}},
$$
and by \eqref{powerreal} together with \eqref{stj}
$$
\norm{\id : [\ell_2^{n-k+1},E_{n-k+1}]_\theta \Id E_{n-k+1}}  \le \norm{\id: \ell_2^{n-k+1} \Id
E_{n-k+1}}^{1-\theta} \prec \left ( \frac{\lambda_E(n-k+1)}{(n-k+1)^{1/2}}
\right )^{1-\theta},
$$
which gives the claim. \qed
\begin{cor}
Let $1 <u_1 <v_1<2$ and $1\le u_2 \le v_2 \le 2$ be such that either $u_2=v_2=2$ or
$\frac{1/{v_1}-1/2}{1/{u_1}-1/2} = \frac{1/{v_2}-1/2}{1/{u_2}-1/2}$. Then for $1 \le k \le n$
\begin{equation}
\label{lorentzck}
c_k(\id: \ell_{v_1,v_2}^n \Id
\ell_{u_1,u_2}^n) \asymp a_k(\id: \ell_{v_1,v_2}^n \Id \ell_{u_1,u_2}^n)
 \asymp (n-k+1)^{1/{u_1}-1/{v_1}}.
\end{equation}
Moreover, formula \eqref{lorentzck} also holds in the case $1 <u_1<v_1<2$ and
$1 \le u_2 \le 2 \le v_2 \le \infty$.
\end{cor}
\proof The first part is clear (use the preceding proposition together with what was mentioned
in the proof of Corollary~\ref{lorentzbc}).
The lower estimates for the second part now follow by factorization:
$$
c_k(\id: \ell_{v_1,v_2}^n \Id \ell_{u_1,u_2}^n) \ge c_k(\id : \ell_{v_1,2}^n \Id \ell_{u_1,2}^n).
$$
The upper estimates are
again straightforward by real interpolation: Choose $0<\theta<1$ such that $1/{v_1}=(1-\theta)/2
+\theta/{u_1}$, then (with the help of \eqref{powerreal})
\begin{samepage}
\begin{align*}
c_k(\id: \ell_{v_1,v_2}^n \Id
\ell_{u_1,u_2}^n) &\le a_k(\id: \ell_{v_1,v_2}^n \Id \ell_{u_1,u_2}^n) \\
&\le \norm{\id : \ell_{v_1,v_2}^{n-k+1} \Id \ell_{u_1,u_2}^{n-k+1}} \\
& \asymp \norm{\id: (\ell_2^{n-k+1},\ell_{u_1,u_2}^{n-k+1})_{\theta,v_2} \Id \ell_{u_1,u_2}^{n-k+1}}
\\ & \le \norm{\id:\ell_2^{n-k+1} \Id \ell_{u_1,u_2}^{n-k+1}}^{1-\theta} \\
& \prec (n-k+1)^{1/{u_1}-1/{v_1}}.
\end{align*}
\qed
\end{samepage} \\[10pt]
We conjecture that formula \eqref{lorentzck} is true for all $1 <u_1 <v_1<2$ and
 $1 \le u_2 ,v_2 \le \infty$.
\providecommand{\bysame}{\leavevmode\hbox to3em{\hrulefill}\thinspace}

Address of the first and the third named author:
\\[10pt]
Fachbereich Mathematik \\ Carl von Ossietzky University of Oldenburg \\ Postfach 2503 \\ D-26111 Oldenburg \\ Germany
\\[5pt]
{\tt defant}@{\tt mathematik.uni-oldenburg.de} \\ {\tt michels}@{\tt mathematik.uni-oldenburg.de}
\\[10pt]
Address of the second named author:
\\[10pt]
Faculty of Mathematics and Computer Science \\ A. Mickiewicz University \\ Matejki 48/49 \\ 60-769 Pozna\'n \\ Poland
\\[5pt] and
\\[5pt]
\begin{samepage}
Institute of Mathematics (Pozna\'n branch) \\ Polish Academy of Sciences \\ Matejki 48/49 \\ 60-769 Pozna\'n \\[5pt]
{\tt mastylo}@{\tt amu.edu.pl}
\end{samepage}

\end{document}